\documentclass[12pt]{article}

\usepackage{amsmath,amssymb,amsbsy,upref}
\usepackage{amsthm}
\usepackage[T2A]{fontenc}
\usepackage[cp1251]{inputenc}
\usepackage{graphicx}
\usepackage{bm}

\newtheorem{Theorem}{Theorem}[section]

\newtheorem{Lemma}{Lemma}[section]

\newtheorem{Remark}{Remark}[section]
\newcommand{\banm}{\begin{anm}}
\newcommand{\eanm}{\end{anm}}

\def\e{\varepsilon}
\def\k{\kappa}

\begin{document}
\title{A fractal graph model of  capillary type systems}
\author{Vladimir Kozlov$^a$, Sergei Nazarov$^b$ and German Zavorokhin$^c$}
\date{}
\maketitle

\vspace{-6mm} 

\begin{center}
$^a${\it Department of Mathematics, Link\"oping University, \\ S--581 83 Link\"oping, Sweden E-mail: vlkoz@mai.liu.se\\
$^b$ St.Petersburg State University,
Universitetsky pr., 28, Peterhof, St. Petersburg, 198504, Russia \\
St. Petersburg State Polytechnical University
Polytechnicheskaya ul., 29, St. Petersburg, 195251, St.Petersburg, Russia \\
Institute of Problems of Mechanical Engineering RAS\\
laboratory “Mathematical Methofs in Mechanics of Materials”\\
V.O., Bolshoj pr., 61, St. Petersburg, 199178, Russia

$^c$ St.Petersburg Department of the Steklov Mathematical Institute, Fontanka, 27, 191023, St.Petersburg, Russia E-mail: zavorokhin@pdmi.ras.ru}

\end{center}


\bigskip \noindent {\bf Abstract.} We consider blood flow in a vessel with an
attached capillary system. The latter is modeled with the help of a corresponding
fractal graph whose edges are supplied with ordinary differential equations obtained
by the dimension-reduction procedure from a three-dimensional model of blood flow in
thin vessels. The Kirchhoff transmission conditions must be satisfied at each
interior vertex. The geometry and physical parameters of this system are described
by a finite number of scaling factors which allow the system to have
self-reproducing solutions. Namely, these solutions are determined by the factors'
values on a certain fragment of the fractal graph and are extended to its rest part
by virtue of these scaling factors. The main result is the existence and uniqueness
of self-reproducing solutions, whose dependence on the scaling factors of the
fractal graph is also studied. As a corollary we obtain a relation between the
pressure and flux at the junction, where the capillary system is attached to the
blood vessel. This relation leads to the Robin boundary condition at the junction
and this condition allows us to solve the problem for the flow in the blood vessel
without solving it for the attached capillary system.

{\it Keywords and phrases}: fractal graph, blood vessel, capillary system, percolation, quiet flow, ideal liquid, Reynolds equation.

\section{Introduction}
\begin{figure}[ht!]
\center{\includegraphics[scale=0.70]{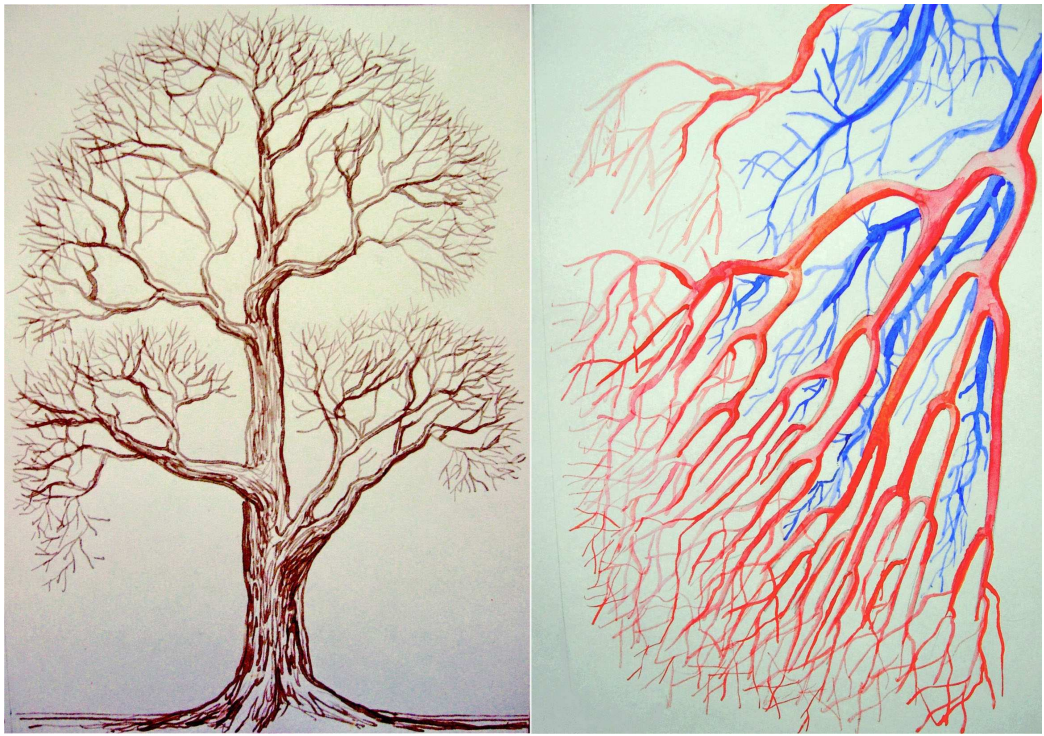}}\\
\caption{A tree and capillaries}\label{Oak}
\end{figure}

{\it Fractal structures} are often involved into both, natural and artificial,
objects in which their elements are repeated iteratively in one or more directions
and simultaneously scaled. In this paper, we study {\it fractal graphs} modelling
capillary blood systems in animal and human bodies as well as vegetative systems in
land-plants and their leaves (see Fig.\ref{Oak}). 

We begin from considering a capillary system as a three-dimensional system of bifurcating thin vessels attached to a main blood vessel at a certain input location. We admit a certain flux through the boundary of the thin vessels and assume that the lengths and the radiuses of the vessels become smaller when we
move away from the main vessel. This will be described with the help of two sets of scaling factors $l_1,\ldots,l_J$ and $k_1,\ldots,k_J$ for the lengths and cross sections  respectively. Passing to the limit, when the radiuses of cross sections tends to zero, we arrive to a one dimensional model (see Appendix for a formal limit procedure). Let us describe this one dimensional model which will serve as our fractal model for the capillary system.

\begin{figure}[ht!]
\center{\includegraphics[scale=0.49]{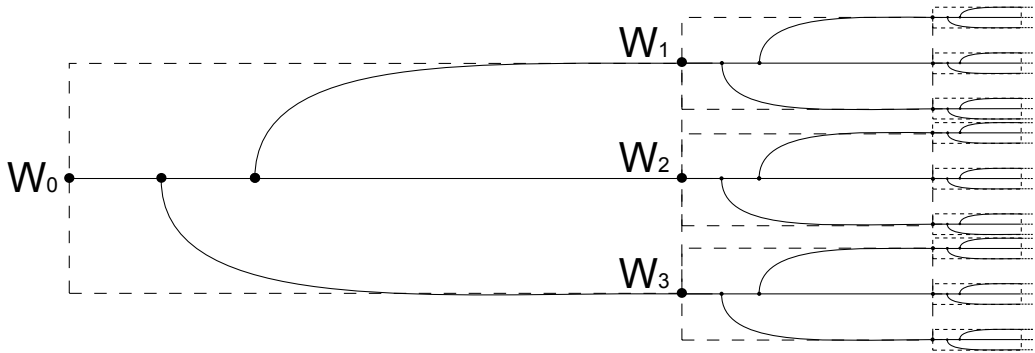}}\\
\caption{Fractal graph $G$ (the case of $J=3$ outputs)}\label{G}
\end{figure}

\begin{figure}[ht!]
\center{\includegraphics[scale=0.50]{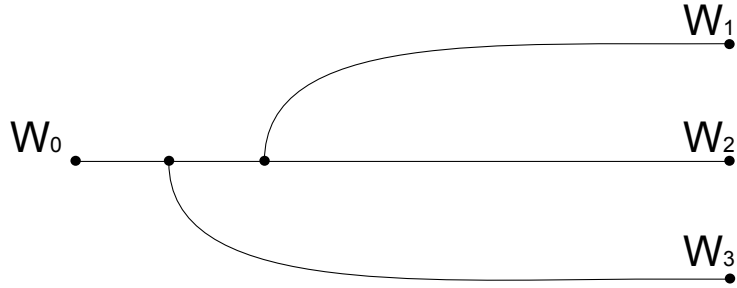}}\\
\caption{Elementary cell $G^0$ (the case of $J=3$ outputs)}\label{G^0}
\end{figure}

Our model includes a fractal graph $G$ (see Fig.\ref{G}). It is obtained from a connected graph $G_0$ with finite number of vertexes and edges
(see Fig.\ref{G^0}). The boundary vertexes (serving as ends only for one edge) consist of one input $W_0$ and $J$ outputs $W_1,\ldots, W_J$.
This graph is attached to the main blood vessel at $W_0$. Our fractal graph is defined by $G_0$ and a given set of scaling factors $l_1,\ldots,l_J\in (0,1)$ in the following way. We attach the graph $l_kG_0$ (with input denoted by $W_{k0}$ and outputs $W_{kj}$, $j=1,\ldots,J$)  to the output $W_k$ of the graph $G_0$ by the input $W_{k0}$. Now we have a  graph with the input $W_0$ and outputs $W_{kj}$. Now we can attach the graphs $l_kl_jG_0$ (with the input denoted by $W_{kj0}$) to the output $W_{kj}$
by the input $W_{kj0}$. Continuing this procedure we obtain our fractal graph, see Fig.\ref{G}.

Next we supply each edge with a differential equation by starting from the graph $G_0$ and extending these differential equations on the other edges of the fractal graph by using scaling factors $l_j$ and  $k_j$, $j=1,\ldots,J$, (see Sect. \ref{Sec.2.2} for more details). Thus
to every edge $e$ of the graph $G$
corresponds a differential equation
\begin{equation}\label{01}
-\partial_{\xi}\left(H(\xi)\partial_{\xi}w(\xi)\right)+ B(\xi)w(\xi)=0,\quad \xi \in e,
\end{equation}
where $\xi$ is the arc length along the edge $e$.
At the vertices the continuity conditions are imposed along with the classical
Kirchhoff transmission condition. This problem on a graph serves as a
one-dimensional model that describes flow of a fluid
in a three-dimensional system of thin channels with the limiting geometry $G$
supplied with boundary conditions on the lateral surface that describes fluid
percolation through the wall of the channel, see Appendix.
In equation (\ref{01}), the unknown $w$ is the hydrodynamic pressure distribution
along a thin channel axis, the given real-valued functions $H(\xi)>0$ and $B(\xi)$
are smooth and describe the throughput capacity of the inferred cross-section
$\omega(\xi)$ of the channel and the total flux through the wall $\partial
\omega(\xi)$ at the point $\xi \in e$, respectively. The Reynolds type equation
(\ref{01}) is a suitable model for steady flow in a thin pipe. This model works for
both, an ideal liquid (the Neumann problem for the three-dimensional Laplace
equation) and a viscous incompressible fluid (the spatial Navier-Stokes equations
with the Robin boundary condition). In the first case, the coefficient $H(\xi)$ is
proportional to the cross-section area, whereas it is proportional to the torsional
rigidity of $\omega(\xi)$ in the second case, see Appendix.

The coefficient $B(\xi)$ is related to either the outgoing (when negative) or
incoming (when positive) flux through the channel's wall. The physical meaning of
both values, negative and positive, is as follows. If a vital wall serves to lead
blood or succus out of the vessel, the outflow through $\partial \omega(\xi)$ is
proportional to the pressure $w(\xi)$ at the point $\xi$, and this gives the term
$B(\xi)w(\xi)$, where $B(\xi)<0$.
In contrast, if an abiotic wall
is porous or damaged with microcracks, the interior pressure $w(\xi)$ increases
permeability of the wall, and so, for a saturated surrounding medium, the total
input $B(\xi)w(\xi)$ with $B(\xi)>0$ must be taken into account in equation
(\ref{01}). We will always assume a certain positivity of the quadratic differential form corresponding to (\ref{01}) on $G_0$ see (\ref{K1a}).
This will provide a certain restriction on the class of the coefficients $B$, when $B$ is negative.

\begin{figure}[ht!]
\center{\includegraphics[scale=0.50]{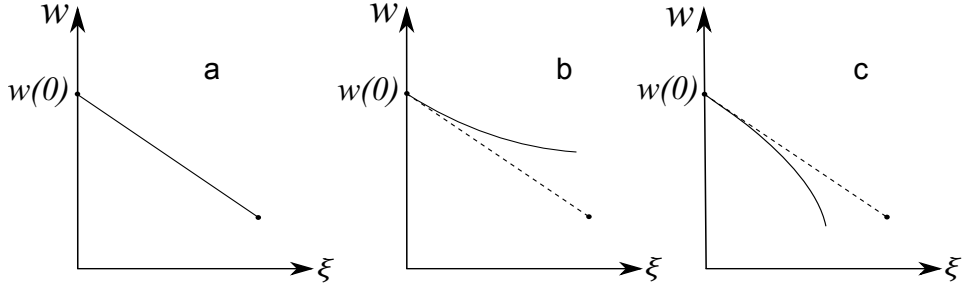}}\\ \caption{Dependence
of pressure $w$ on $\xi$ for $B$ values of different sign: a - flux
is constant $(B=0)$, b - flux decreases $(B>0)$, c - flux increases
$(B<0)$}\label{w(0)}
\end{figure}

A construction
for the differential equations on the edges fractal graph $G$ is
given in Sect.\ref{Sec.2.2}. Here we present just the first step in this
construction.
Let  $k_{j}\in(0,1), j=1,\ldots,J,$ be the scaling factors for
 ``cross-sections''  in the three-dimensional model, see
Appendix for details. We start with the graph $G^0$ and attach the input of the
graph $l_j G^0$ to the $j$-th output of $G^0$. The coefficients $H^{j}_{e}$ and
$B^{j}_{e}$ are defined by
\begin{equation}\label{HBj}
H^{j}_{e}(\xi)=k_{j}H_{e}\left(\frac{\xi}{l_{j}}\right),\
B^{j}_{e}(\xi)=\frac{k_{j}}{l^{2}_{j}}B_{e}\left(\frac{\xi}{l_{j}}\right),\ \xi \in
[0,L_{e}l_{j}]
\end{equation}
for every edge $e$ of the graph $l_j G^0$. It is clear that if $w$ is a solution to
(\ref{01}) with the Kirchhoff transmission conditions which, in particular, include
the continuity at the interior vertices of $G^0$, then
\begin{equation}\label{002}
w^j(\xi)=m_{j}w\left(\frac{\xi}{l_{j}}\right)
\end{equation}
is a solution of the corresponding problem on $l_j G^0$, where the constants $m_j$
are arbitrary. In this way, only self-reproducing solutions are sought, i.e.
solutions of the form (\ref{002}) which are continuous at the junction points, where
the Kirchhoff transmission condition must be satisfied above all. This leads to the
following equations for $m_j$:
\begin{equation}\label{003}
w(W_j)=m_jw(W_0),
\end{equation}
\begin{equation}\label{004}
F_j(w)=m_j\kappa_jF_0(w),
\end{equation}
where $j=1,\ldots,J, \kappa_j=\frac{k_j}{l_j}$, and $$F_j(w)=-H(W_j)w'(W_j)$$ is the
flux at the vertex $W_j$; the outward derivative with respect to the graph $G^0$ is
taken at $F_1,\ldots,F_J$, whereas it is directed inwards at $F_0$. If we put
$F_j(X_0,{\bf X})=F_j(w)$, where ${\bf X}=(X_1,\ldots,X_J)$ and $w$
satisfies the Dirichlet boundary conditions $w(W_0)=X_0$ and $w(W_j)=X_j$,
$j=1,\ldots,J$, then systems (\ref{003}) and (\ref{004}) can be written as the
following single system:
\begin{equation}\label{004a}
F_j(1,{\bf m})=m_j\kappa_jF_0(1,{\bf m}),\;\;\;j=1,\ldots,J.
\end{equation}
These relations are essential for finding ${\bf m}$. When ${\bf m}$ is found we
obtain the following relation
\begin{equation}\label{004b}
F_0(w)=\beta({\bf m})w(W_0),\;\;\;\beta({\bf m})=F_0(1,{\bf m})
\end{equation}
between the Dirichlet and Neumann boundary conditions at the input $W_0$ for
arbitrary self-reproducing solution with parameters $k_j$, $l_j$.

Most of computational schemes are difficult to implement for graphs of fractal types
that have a complicated topological arrangement. In particular, this happens because
it is necessary to reduce permanently the branching tree of the grid spacing or the
size of the spline. However, many  objects containing fractal fragments have common
features:
  veins and arteries are involved in the movement of blood on a large scale to and from various parts of the body, whereas capillaries are involved in the local distribution of blood to cells and body's tissues on a small scale. Similar versions of liquid distribution occur in vegetative systems of plants and their fragments.

Our goal is to replace the fractal branch attached to a "main" vessel by artificial
boundary conditions (\ref{004b}) imposed where the attachment is localised. To this
end, the main problem is to determine the parameter $\beta({\bf m})$ or,
equivalently, ${\bf m}$. It appears that this parameter can be found under
reasonable assumptions on the fractal graph.

In what follows, we look for a positive solution to the problem described above,
that is, all $m_j>0$. The case of $m_j<0$ corresponds to extraction of fluid out
of the system, but this does not occur for vegetative and capillary systems. 
We also assume that $m_j<1$, because the system under consideration will have an
unbounded growth of pressure in the channels otherwise. This is in conflict with the
normal system's functioning and its viability.

The amount of small capillaries in human body is about several milliard, and so it
is hardly possible to represent them as a one-dimensional network with a
differential structure. Besides, there are two main threats to functioning of the
entire circulatory system properly, namely: clogging and damages of the arterial
tree (blood clots, aneurysms, crushing vessel walls, venous nodes etc.). However,
even disruption of the venous valves or bleeding from veins (gemoroin phenomena) are
not immediately fatal, not to mention the injury of large groups of capillaries
(external bruising, internal hematoma etc) the latter can happen dayly. Therefore,
it is not apparently necessary to solve explicitly the problem of blood flow in the
capillary system which describes distribution of blood in human body. Nevertheless,
it is important to model arterial blood outflow due to capillaries.

In the current literature, there are several approaches to modeling the capillary
network or its parts. In particular, the paper \cite{Poz} deals with numerical
modeling of branching microvascular tree-like capillary network, where the linear
Poiseuille flow is adopted (see the Reynolds equation (\ref{A15}) and (\ref{strong})
below). However, the effective viscosity of blood depends on the corrsponding
distribution of red cells (the dependence is found out from other numerical
experiment on individual red blood cell motion). In the paper \cite{Dav}, the linear
stability of this model is studied in the tree-like and honey comb networks, for
which purpose many numerical experiments are presented that describe miscellaneous
particular effects of blood flow in the branching capillary systems. A different,
but to some extent similar model is proposed in \cite{FiLaLi}; it takes into account
the elastic properties of the walls of capillaries, external influence on their
surrounding muscle and blood seeping through the walls.

The closest to the spirit of our work is the paper \cite{APQ}, where averaging
method is applied for studying what is brought into the one-dimensional Reynolds
equation by the flow of blood from the main vessel to a system of thin vessels that
are parallel to each other and arranged periodically. We also use a version of
averaging method, but consider the capillary network as a fractal tree which is
aimed to finding a connection between the blood pressure and the loss of flow in
the capillary tree.

Also, a lot of works deals with flows on networks; see, for example, the survey paper
\cite{BCG}. The novelty of our study is to modeling the capillary system as a
fractal network described by a set of scaling factors and the corresponding
self-reproducing solutions for blood flow in such a fractal structure.

Description of the fractal graph $G$ and statement of the problem are given in Sect.
\ref{Sec.2}. Sect. \ref{S3} concerns the study of equation (\ref{004a}). Our
analysis is splitted between three cases: (i) impermeable case ($B=0$); (ii)
permeable case ($B\geq 0$); (iii) permeable case ($B\leq 0$). The study is based on
the maximum principle. In the case (i), we show that equation (\ref{004a}) is
uniquely solvable for all $k_j$, $l_j$ such that $k_1l_1^{-1}+\cdots+
k_Jl_J^{-1}>1$. In the case (ii), the unique solvability is proved for all positive
$k_j$ and $l_j$. Finally, in the case (iii), which actually occurs in capillary
systems, we prove that (\ref{004a}) has a unique solution provided its energy is
positive, i.e.
 $$
 F_0(1,{\bf m})-\sum_{j=1}^JF_j(1,{\bf m})m_j>0.
 $$
Without the last assumption, there is no uniqueness. A derivation of one-dimensional
models from three-dimensional ones is given in Appendix.

\section{Formulation of the problem}\label{Sec.2}
\subsection{Model problem on elementary cell $G^0$}\label{Sec. 2.1}

Let $(V,E)$ be a connected graph in $\mathbb R^3$ with vertexes $V$ and edges $E$. We denote by $W=\{W_j\,:\,j=0,1,\ldots,J\}, J\geq 1$, the vertexes which are attached only to one edge. All others are attached at least to two edges. We call the vertex $W_0$ input and the vertexes $W_j$, $j=1,\ldots,J$, outputs. We represent each edge $e\in E$ as a curve with corresponding vertexes as endpoints and we supply it with the natural parametrization $\xi\in (0,L_e)$, where $L_e$ is the length of the curve. It is supposed that  functions $H_e\in C^1[0,L_e]$ and $B_e\in C[0,L_e]$ are given on each edge $e\in E$. It is assumed that there exists a positive constant $c_H$  such that
$$
c_H\leq \min_{\xi\in [0,L_e]} H_e(\xi)\;\;\;\mbox{for all $e\in E$}.
$$

Let $\overline{G^0}$ be the union of all edges and vertexes in the graph $(V,E)$ and let $G^0=\overline{G^0}\setminus\{W_0,W_1,\ldots,W_J\}$.
By $H^1(G^0)$ we denote  the set of real-valued functions $w$ on $G^0$, continuous on each edge and at the vertexes with the norm
$$
||w||_{H^1(G^0)}=\Big (\sum_{e\in E}\int_0^{L_e}(|w'(\xi)|^2+|w(\xi)|^2)d\xi\Big)^{1/2},
$$
where $w'$ is the derivative of $w$ with respect to $\xi$.
Let also $H^1_0(G^0)$ be the subspace in $H^1(G^0)$ consisting of functions equal zero at all points $W_j$, $j=0,1,\ldots,J$.
In order to define a weak formulation of the problem which we are going to study, we introduce a bilinear form on $H^1(G^0)$:
\begin{equation}\label{00.1}
a(w,u)=\sum_{e\in E}\int_0^{L_e}(Hw'(\xi)u'(\xi)+Bw(\xi)\,u(\xi))d\xi.
\end{equation}
Then we introduce the following Dirichlet problem: find $w\in H^1(G^0)$ such that
\begin{equation}\label{00.2}
a(w,u)=0\;\;\;\mbox{for all $u\in H^1_0(G^0)$},
\end{equation}
satisfying the Dirichlet boundary conditions
\begin{equation}\label{00.3}
w(W_j)=X_j,\;\;\;j=0,1,\ldots,J.
\end{equation}
One can use an equivalent strong formulation of problem (\ref{00.2}), (\ref{00.3}):
\begin{equation}\label{strong}
-(Hw')'+B w=0\;\;\mbox{on $e$},
\end{equation}
$w$ is continuous at each vertex, the Kirchhoff conditions are valid at each vertex of $V\setminus W$ and (\ref{00.3}) is satisfied. We recall that  the Kirchhoff condition at $v\in V\setminus W$ is defined as
$$
\sum H_e(v)w'(v)=0,
$$
where the sum is taken over all $e$ attached to the vertex $v$ and the derivative is taken outwards with respect to the vertex $v$.


We always assume that
\begin{equation}\label{K1a}
a(u,u)>0\;\;\mbox{for all $u\in H^1_0(G^0)\setminus \{O\}$.}
\end{equation}
Then the problem (\ref{00.2}), (\ref{00.3})
has a unique solution. As a consequence, let us prove the following
\begin{Lemma}\label{QuQu1} Let $X_j\geq 0$, $j=0,\ldots,J$, and let at least one of them be positive.
Also let $w$ be a solution of the Dirichlet problem {\rm (\ref{00.2})}, {\rm (\ref{00.3})}.
Then $w>0$ in $G^0$. If additionally $X_k=0$ for certain $k=0,\ldots,J$, then $F_k(X_0,{\bf X})<0$ for $k=0$ and $F_k(X_0,{\bf X})>0$ for $k>0$.
\end{Lemma}
{\bf Proof.} Let us show that this solution $w$ is positive.
We note that the form $a=a_t$, which is obtained from $a$ if $B$ is replaced by $tB$,
satisfies also (\ref{K1a}) when $t\in [0,1]$. Therefore, one can define also the solution $w_t$ to (\ref{00.2}) subject to the boundary conditions (\ref{K1b}). One can readily check that $w_0>0$
in $G^0$ by the maximum principle and that $w_t$ continuously depends on $t\in[0,1]$. Let us take the first $t$ (we denote it by $t_0$) for which $w_{t_0}$ has  zero at a certain interior point. Since this is the minimum, then  the derivative also vanishes at the same point. Due to uniqueness of solutions to the Cauchy problem the function vanishes on the whole edge containing this point together with vertexes-endpoints of the edge. Since zero is the minimum of the function then using the Kirchhoff boundary condition one can show that all fluxes at these vertexes are zero. Therefore, the function $w_{t_0}$ vanishes at edges adjacent to these vertexes. Continuing this procedure, we obtain that $w_{t_0}=0$ on $G^0$. This contradiction ($w$ is not identically zero by the assumptions of the lemma) demonstrates that $w_t$ is a positive function in $G^0$ for all $t\in [0,1]$ and, in particular, for $t=1$.

To prove the assertions about the flux, we can use the same argument to show that the flux cannot be zero at $W_k$ for all $t\in [0,1]$.

\medskip
The following Green formula will be used in the sequel:

\begin{equation}\label{LUD1}
a(w,v)=F_0(X_0,{\bf X})Y_0-\sum_{j=1}^JF_j(X_0,{\bf X})Y_j,
\end{equation}
where $w,v\in H^1(G^0)$ are solutions to (\ref{00.2}) with the Dirichlet
data $(X_0,{\bf X})$ and $(Y_0,{\bf Y})$ respectively. As a
consequence, we get
\begin{equation}\label{LUD1b}
\sum_{j=1}^JF_j(X_0,{\bf X})Y_j-F_0(X_0,{\bf X})Y_0=
\sum_{j=1}^JF_j(Y_0,{\bf Y})X_j-F_0(Y_0,{\bf Y})X_0.
\end{equation}

We will distinguish three cases: a) $B\equiv 0$ (impermeable walls); b) $B\geq 0$ and $B$ is not identically zero (permeable walls); c) $B\leq 0$ and $B$ is not identically zero (permeable walls also).

\subsection{Fractal graph $G$}\label{Sec.2.2}

In order to describe the fractal graph $G$ we introduce numbers $l_j \in (0,1), j=1,\ldots,J,$ and put
$$
G(J_n)=l_{j_1}l_{j_2}\ldots l_{j_n}G^0,
$$
where $J_n=(j_1,\ldots,j_n)$, and $j_k$ may take values $1,\ldots,J$.
We denote by $W_{0}^{J_n}$ the input of  $G(J_n)$ and by $W^{J_n}_{j}$ its outputs. Further, we identify the input $W_{0}^{J_n}$ of the graph $G(J_n)$ with the output $W_{j}^{J_{n-1}}$ of the graph $G(J_{n-1}), J_{n-1}=(j_1,\ldots,j_{n-1}),$ (we put $G(J_0)\equiv G^0$). These graphs together with above identification  give the fractal graph
$$
G=\bigcup_{n=0}^{+\infty}G(J_n).
$$
To define the corresponding differential operators we introduce numbers $k_j \in (0,1), j=1,\ldots,J,$ and on each edge $e$ of the graph $G(J_n)$ we define the functions
$$
H^{J_n}_{e}(\xi)=k_{j_1}k_{j_2}\ldots k_{j_n}H_{e}\left(\frac{\xi}{l_{j_1}l_{j_2}\ldots l_{j_n}}\right),\
$$
$$
B^{J_n}_{e}(\xi)=\frac{k_{j_1}k_{j_2}\ldots k_{j_n}}{l^{2}_{j_1}l^{2}_{j_2}\ldots l^{2}_{j_n}}B_{e}\left(\frac{\xi}{l_{j_1}l_{j_2}\ldots l_{j_n}}\right),\
$$
where $\xi \in [0,Ll_{j_1}l_{j_2}\ldots l_{j_n}].$

Note that fractal graph $G$ consists of edges that are closed arcs. It can have cycles and occupies finite volume, since due to scaling factors $l_{j} \in (0,1)$ the sizes of cells decrease as a geometrical progression.

We are looking for a solution to the model problem (\ref{00.2}), (\ref{00.3}) on $G^0$ which has self-reproducing structure on the whole fractal graph $G$.

Namely, if $w(\xi)$ is a solution to the problem (\ref{00.2}), (\ref{00.3}) on $G^0$ then
$$
w^{J_n}(\xi)=m_{j_{1}}m_{j_{2}}\ldots m_{j_{n}}w\left(\frac{\xi}{l_{j_{1}}l_{j_{2}}\ldots l_{j_{n}}}\right)
$$
is a solution of the corresponding problem on the graph $G(J_n)$ for arbitrary coefficients $m_j$. This self-reproducing solution solves the problem on the whole graph $G$ if the continuity and the Kirchhoff transmission conditions are satisfied at all attachment points. Thus we obtain  equations (\ref{004a}) for the scaling factors $m_j$.
Moreover, the vector ${\bf m}$ is sought in the set
\begin{equation}\label{K3}
\Omega\!\!=\!\!\left\{{\bf m}\in \mathbb{R}^{J}\;|\;F_{0}(1,{\bf m})>0,
\; 1>m_{j}>0,  F_{j}(1,{\bf m})> 0, j=1,\ldots,J\right\}\!.
\end{equation}
Determination of $m_j$ gives a possibility to obtain the relation between the Dirichlet and Neumann data at the input $W_0$, see (\ref{004b}).
This relation represents the Robin boundary condition with the coefficient $\beta$.

In what follows, our main concern is to study solvability of equation (\ref{004a}), in particular, to describe  the set of $\kappa_j$ for which there exist solutions ${\bf m}$ and to study their multiplicity.

\section{Auxiliary assertions}

Introduce the vector function ${\cal F}=({\cal F}_1,\ldots,{\cal F}_J)$, where
\begin{equation}\label{k1}
{\cal
F}_{j}={\cal F}_j({\bf m})=\frac{F_{j}(1,{\bf m})}{F_{0}(1,{\bf m})m_{j}},
\end{equation}
which is considered on $\Omega$ or
\begin{equation}\label{K2}
\widehat{\Omega}\!\!=\!\!\left\{{\bf m}\in \mathbb{R}^{J}\;|\;F_{0}(1,{\bf
m})>0 ,\; 1>m_{j}>0,  F_{j}(1,{\bf m})\geq 0.
j=1,\dots,J\right\}.
\end{equation}
Since these sets are given by linear inequalities both of them are
convex and $\widehat{\Omega}\supset\Omega$. Now equations (\ref{004a}) can be written as
\begin{equation}\label{03newQ}
{\cal F}_j({\bf m})=\kappa_j, \;j=1,\ldots,J.
\end{equation}
Our aim is to investigate solvability of this system of equations.

The Jacobian of the vector function ${\cal F}$ is denoted by
$$
{\cal J}({\bf m})=\{{\cal J}_{ji}({\bf m})\}_{j,i=1}^J,
$$
where
$$
{\cal J}_{ji}({\bf m})=\frac{F_j(0,{\bf e}_i)-{\cal F}_j({\bf m})F_0(0,{\bf e}_i)m_j-{\cal F}_j({\bf m})
F_0(1,{\bf m})\delta_i^j}{F_0(1,{\bf m})m_j}.
$$
Here $\delta_i^j$ is the Kronecker delta and ${\bf
e}_k=\{\delta_k^j\}_{k=1}^J$.
We introduce the following $(J-1)$-dimensional subspace of $\mathbb R^J$
\begin{equation}\label{TTa}
{\cal R}_0=\{ {\bf H}\in\mathbb R^J\,:\, F_0(0,{\bf H})=0\}.
\end{equation}
In the next lemma we give a necessary and sufficient condition for invertibility of the Jacobian ${\cal J}$.

\begin{Lemma}\label{Lem1} Let us assume that for all ${\bf m}\in\widehat{\Omega}$
\begin{equation}\label{TTf1s}
a(V,V)+\sum_{j=1}^JF_j(1,{\bf m})m_j^{-1}H_j^2>0\;\;
\end{equation}
for all  ${\bf H}\in {\cal R}_0\setminus\{ O\}$. Here $V$ is the solution to   {\rm (\ref{00.2})}, {\rm (\ref{00.3})}
with $(X_0,{\bf X})=(0,{\bf H})$.
Then the following relation holds:
\begin{equation}\label{TTf1}
\det{\cal J}({\bf m})=\Lambda({\bf m})\Big(F_0(1,{\bf m})-\sum_{j=1}^JF_j(1,{\bf m})m_j\Big),\;\;{\bf m}\in\widehat{\Omega},
\end{equation}
where $\Lambda$ is smooth, non-vanishing function on $\widehat{\Omega}$.
\end{Lemma}
{\bf Proof.} First, we note that
$$
\det{\cal J}=\frac{\det{\cal S}}{F_0(1,{\bf m})^Jm_1\cdots m_J},
$$
where ${\cal S}=\{ {\cal S}_{ji} \}_{j,i=1}^J$ and
$$
{\cal S}_{ji}=F_j(0,{\bf e}_i)-{\cal F}_j({\bf m})F_0(0,{\bf e}_i)m_j-{\cal F}_j({\bf m})
F_0(1,{\bf m})\delta_i^j.
$$
Furthermore, if ${\bf h},\,{\bf g}\in\mathbb R^J$, then
\begin{eqnarray*}
&&\sum_{1\leq i,j\leq J}{\cal S}_{ji}h_ig_j=\sum_{j=1}^J\big(F_j(0,{\bf h})g_j-F_0(1,{\bf m}){\cal F}_jh_jg_j-F_0(0,{\bf h}){\cal F}_jm_jg_j\big)\\
&&=\sum_{j=1}^J\big((F_j(0,{\bf g})-F_0(1,{\bf m}){\cal F}_jg_j)h_j-F_0(0,{\bf h}){\cal F}_jm_jg_j\big),
\end{eqnarray*}
where we have used the relation $\sum F_j(0,{\bf h})g_j=\sum F_j(0,{\bf g})h_j$, which follows from (\ref{LUD1b}). If we take ${\bf g}={\bf m}$, we get
\begin{equation}\label{UU1a}
\!\!\sum_{1\leq i,j\leq J}\!\!{\cal S}_{ji}h_im_j\!=\!-\sum_{j=1}^J(F_j(1,{\bf 0})h_j+F_0(0,{\bf h}){\cal F}_jm_j^2)=F_0(0,{\bf h})\Big(1-\sum_{j=1}^J{\cal F}_jm_j^2\Big),
\end{equation}
where the relation $\sum F_j(1,{\bf 0})h_j=F_0(1,{\bf 0})$ is used, which again follows from (\ref{LUD1b}).

Similarly, if ${\bf h}={\bf m}$ then
\begin{equation}\label{UU1b}
\sum_{1\leq i,j\leq J}{\cal S}_{ji}m _ig_j=F_0(0,{\bf g})-F_0(0,{\bf m}){\cal F}_jm_jg_j.
\end{equation}
By Lemma \ref{QuQu1}, $F_0(0,{\bf m})<0$ for ${\bf m}\in\widehat{\Omega}$. Therefore, each vector in $\mathbb R^J$ can be represented as a linear combination of a vector in ${\cal R}_0$ and ${\bf m}$ and this representation is unique. Let
$$
{\bf h}={\bf H}+\lambda{\bf m},\;\;{\bf g}={\bf G}+\mu{\bf m},\;\;{\bf H},\,{\bf G}\in {\cal R}_0\;\;\mbox{and}\;\;\lambda,\,\mu\in\mathbb R.
$$
Using (\ref{UU1a}) and (\ref{UU1b}), we get
\begin{eqnarray}\label{UU1c}
&&\sum_{1\leq i,j\leq J}{\cal S}_{ji}h_ig_j=\sum_{1\leq i,j\leq J}{\cal S}_{ji}H_iG_j+\lambda\mu F_0(0,{\bf m})\Big(1-\sum_{j=1}^J{\cal F}_jm_j^2\Big)\nonumber\\
&&-\lambda F_0(0,{\bf m})\sum_{j=1}^J{\cal F}_jm_jG_j.
\end{eqnarray}
Let us show that the bilinear form
$$
b({\bf H},{\bf G})=\sum_{1\leq i,j\leq J}{\cal S}_{ji}H_iG_j,\;\;{\bf H},{\bf G}\in {\cal R}_0,
$$
is non-degenerate. Indeed, assume that there exist nontrivial  ${\bf H}$ such that $b({\bf H},{\bf G})=0$ for all ${\bf G}\in {\cal R}_0$. Since
$$
b({\bf H},{\bf G})=\sum_{j=1}^J\big(F_j(0,{\bf H})-F_0(1,{\bf m}){\cal F}_jH_j\big)G_j\;\;\mbox{and}\;\;\sum_{j=1}^JF_j(1,{\bf 0})G_j=-F_0(0,{\bf G})=0,
$$
we conclude that
\begin{equation}\label{UU1d}
F_j(0,{\bf H})-F_0(1,{\bf m}){\cal F}_jH_j=aF_j(1,{\bf 0}),\;\;j=1,\ldots ,J,
\end{equation}
for a certain constant $a$. Multiplying relation (\ref{UU1d}) by $m_j$ and summing them up, we obtain
$$
\sum (F_j(0,{\bf m})-F_0(1,{\bf m})m_j{\cal F}_j)H_j=-aF_0(0,{\bf m}),
$$
where the relations
$$
\sum F_j(0,{\bf H})m_j=\sum F_j(0,{\bf m})H_j\;\;\;\mbox{ and}\;\;\; F_j(1,{\bf 0})m_j=-F_0(0,{\bf m})
$$
 are used. Since the left-hand side equals $-\sum F_j(1,{\bf 0})H_j=F_0(0,{\bf H})=0$ and $F_0(0,{\bf m})\neq 0$, we conclude that $a=0$. But then, multiplying (\ref{UU1d}) by $H_j$, summing up and using (\ref{LUD1}), we get
\begin{equation}\label{UU1e}
-a(V,V)=\sum F_0(1,{\bf m}){\cal F}_jH_j^2,
\end{equation}
where $V$ is a solution to (\ref{00.2}), (\ref{00.3}) with the Dirichlet data $(X_0,{\bf X})=(0,{\bf H})$. Due to the positivity condition (\ref{TTf1s}), equation (\ref{UU1e}) implies
${\bf H}=0$. Thus the form $b$ is non-degenerate. Now, we can interpret relation (\ref{UU1c}) as representation of the matrix ${\cal S}$ in the block form in a certain basis in $\mathbb R^J$ with one zero block-matrix. Since the block corresponding to the form $b$ is non-degenerate the corresponding block-matrix has non-zero determinant and hence,
$$
\det{\cal S}=C({\bf m})\Big(1-\sum_{j=1}^J{\cal F}_jm_j^2\Big),
$$
where $C$ is non-vanishing, smooth function on $\widehat{\Omega}$. This leads to (\ref{TTf1}). The proof is complete.

\begin{Remark}\label{Rem11} Let us denote by $M$ a solution to {\rm (\ref{00.2})}, {\rm (\ref{00.3})},
where $(X_0,{\bf X})=(1,{\bf m})$. Then, by {\rm (\ref{LUD1})},
\begin{equation}\label{LUD1s}
F_0(1,{\bf m})-\sum_{j=1}^JF_j(1,{\bf m})m_j=a(M,M).
\end{equation}
Hence the second factor in the right-hand side of {\rm (\ref{TTf1})} can be interpreted as an energy of the solution $M$.
\end{Remark}

In the following lemma we present an injectivity criterion for the mapping ${\cal F}$.
\begin{Lemma}\label{Lem1s} We assume that {\rm (\ref{TTf1s})} is valid on $\widehat{\Omega}$. Let ${\bf m},\,{\bf n}\in\widehat{\Omega}$ satisfy
\begin{equation}\label{LUD1ss}
{\cal F}_j({\bf m})={\cal F}_j({\bf n}),\;\;j=1,\ldots,J.
\end{equation}
If
\begin{equation}\label{LUD1sf}
F_0(1,{\bf m})-\sum_{j=1}^JF_j(1,{\bf m})m_j>0\;\;\mbox{and}\;\;F_0(1,{\bf n})-\sum_{j=1}^JF_j(1,{\bf n})n_j>0
\end{equation}
then ${\bf m}={\bf n}$.
\end{Lemma}
{\bf Proof}. Let $\kappa_j={\cal F}_j({\bf m})$. Since ${\bf m}$ and ${\bf n}$ belong to $\widehat{\Omega}$, the numbers $\kappa_j$ are non-negative.
From (\ref{LUD1ss}) it follows
$$
F_j(1,{\bf m})- F_j(1,{\bf n})=\kappa_j(F_0(1,{\bf m})m_j- F_0(1,{\bf n})n_j),\;\;j=1,\ldots,J.
$$
Using notation ${\bf h}={\bf m}-{\bf n}$, we derive from the last relation
\begin{equation}\label{LUD1sa}
F_j(0,{\bf h})=\kappa_j(F_0(0,{\bf h})m_j+F_0(1,{\bf n})h_j),\;\;j=1,\ldots,J.
\end{equation}
The vector ${\bf h}\in\mathbb R^J$ admits the representation (compare with the proof of Lemma \ref{Lem1}) ${\bf h}={\bf H}+\lambda{\bf m}$ where $F_0(0,{\bf H})=0$ and $\lambda\in\mathbb R$.
Consider two cases.

1). ($\lambda=0$) Then (\ref{LUD1sa}) implies $F_j(0,{\bf H})=\kappa_jF_0(1,{\bf n})H_j$. Multiplying these equalities by $H_j$ and summing them up, we get
$$
-\sum_{j=1}F_j(0,{\bf H})H_j+\sum_{j=1}^JF_j(1,{\bf n})n_j^{-1}H_j^2=0,
$$
which is impossible due to (\ref{TTf1s}).

2). ($\lambda\neq 0$) Then (\ref{LUD1sa}) leads to
$$
F_j(0,{\bf H})-\kappa_jF_0(1,{\bf n})H_j=\lambda\kappa_j\big(F_0(1,{\bf m})m_j+F_0(1,{\bf n})m_j\big)-\lambda F_j(0,{\bf m}).
$$
Since
\begin{eqnarray*}
&&\sum_{j=1}^Jn_j(F_j(0,{\bf H})-\kappa_jF_0(1,{\bf n})H_j)=\sum_{j=1}^J(F_j(0,{\bf n})-\kappa_jF_0(1,{\bf n}))H_j\\
&&=-\sum_{j=1}^JF_j(1,{\bf 0})H_j=F_0(0,{\bf H})=0,
\end{eqnarray*}
we have that
\begin{eqnarray}\label{22dec1}
&&0=\sum_{j=1}^Jn_j\Big(\kappa_j\big(F_0,{\bf m})m_j+F_0(1,{\bf n})m_j\big)- F_j(0,{\bf m})\Big)\\
&&=\sum_{j=1}^JF_0(0,{\bf m})\kappa_jm_jn_j+\sum_{j=1}^JF_j(1,{\bf 0})m_j=F_0(0,{\bf m})\Big(\sum_{j=1}^J\kappa_jm_jn_j-1\Big),\nonumber
\end{eqnarray}
where we have used that $\sum F_j(1,{\bf 0})m_j=F_0(0,{\bf m})$ by (\ref{LUD1b}). From (\ref{LUD1sf}) it follows that
$$
\sum_{j=1}^J\kappa_jm_j^2<1\;\;\mbox{and}\;\;\;\sum_{j=1}^J\kappa_jn_j^2<1.
$$
the latter implies
$$
2\sum_{j=1}^J\kappa_jm_jn_j\leq \sum_{j=1}^J\kappa_jm_j^2+\sum_{j=1}^J\kappa_jn_j^2<2.
$$
Hence the right-hand side of (\ref{22dec1})  does not vanish. This proves the required assertion.

\section{Solution of equation (\protect\ref{03newQ})}\label{S3}

\subsection{The case of impermeable walls}\label{Sect.3.1}

\begin{figure}[ht!]
\center{\includegraphics[scale=0.50]{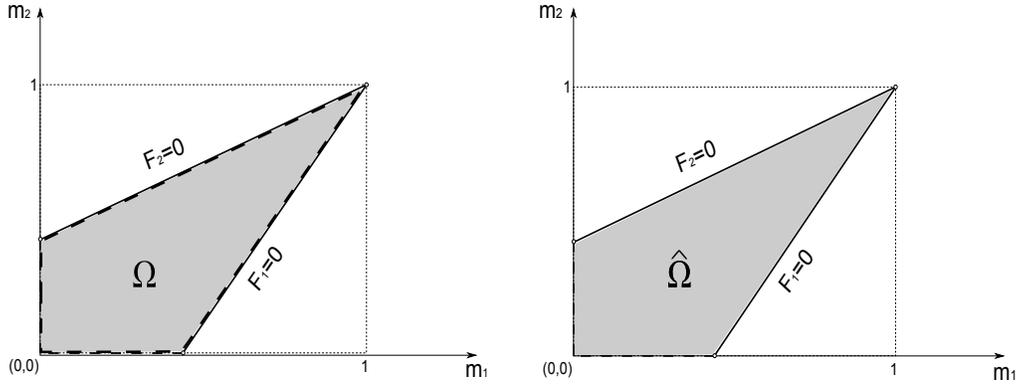}}\\
\caption{Sets $\Omega,\widehat{\Omega}$ (the case of $J=2$
outputs)}\label{omega}
\end{figure}

In this section we assume that  $B\equiv0$.
The sets $\widehat{\Omega}\supset\Omega$, introduced by (\ref{K3}) and (\ref{K2}) respectively,   are outlined on  Fig.\ref{omega}.
Let ${\bf 1}=(1,\ldots,1) \in \mathbb{R}^{J}$. Applying (\ref{LUD1}) to the function $v=1$ corresponding to
the Dirichlet data $(1,{\bf 1})$ and observing that $a(w,v)=0$ for arbitrary $w$ in
this case, we get
\begin{equation}\label{LUD4a}
\sum_{j=1}^JF_j(X_0,{\bf X})=F_0(X_0,{\bf X}).
\end{equation}

In what follows, we will use the following

\smallskip
\noindent {\bf Maximum principle  ($B=0$)}. {\em If $v(\xi)\neq const$ is a
solution to problem (\ref{00.2})-(\ref{00.3}) with the Dirichlet
data $(X_{0},{\bf X})$ , then
$$
\min_{0\leq j\leq J}X_{j}< v< \max_{0\leq j\leq
J}X_{j}\;\;\;\mbox{in the interior of $G^0$.}
$$
Moreover,  if minimum (maximum) is attained  at $W_k$ then $F_k>0$
when $k=1,\ldots,J$ and $F_k<0$ when $k=0$ ($F_k<0$ when
$k=1,\ldots,J$ and $F_k>0$ when $k=0$ in the case of maximum)}.

Let us show that the set $\Omega$ can be described by a less number of inequalities. Namely,
\begin{equation}\label{XX1a}
\Omega=\{ {\bf m}\,:\,m_j>0\;\mbox{and}\; F_j(1,{\bf m})>0,\,j=1,\ldots,J\,\}.
\end{equation}
Indeed, let $M$ satisfy (\ref{00.2})-(\ref{00.3}) with the Dirichlet data $(1,{\bf m})$. Then application of the maximum principle
to the function $1-M$ gives $m_j<1$ (if $m_j>1$ for certain $j$ then $F_j(1,{\bf m})<0$ which is not true for elements in $\Omega$).
From  (\ref{LUD4a}) it follows that $F_0(1,{\bf m})>0$. Similarly,
the set $\widehat{\Omega}$ can be described as
\begin{equation}\label{XX1b}
\widehat{\Omega}=\{ {\bf m}\,:\,m_j>0\;\mbox{and}\; F_j(1,{\bf m})\geq 0,\,j=1,\ldots,J,\,F_0(1,{\bf m})>0\}.
\end{equation}


\begin{Lemma}\label{Lem2} The Jacobian ${\cal J}$ is invertible on $\widehat{\Omega}$.
The map ${\cal F}:\widehat{\Omega}\rightarrow \mathbb{R}^{J}$ is injective.
\end{Lemma}
{\bf Proof.} In the case $B\equiv 0$, $a(v,v)>0$ for nontrivial $v\in H^1(G^0)$ satisfying $v(W_0)=0$, which implies (\ref{TTf1s}).
By Remark \ref{Rem11}, the result follows from Lemmas \ref{Lem1} and \ref{Lem1s}.

\begin{figure}[ht!]
\center{\includegraphics[scale=0.49]{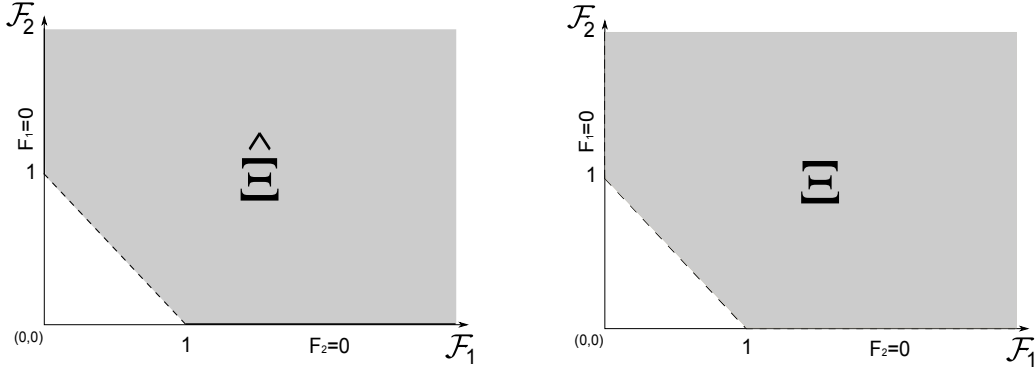}}\\
\caption{Sets $\widehat{\Xi},\Xi$ (the case of $J=2$ outputs)}\label{ksi}
\end{figure}

\medskip
We put (see Fig.\ref{ksi})
$$
\widehat{\Xi}=\left\{{\cal F}\in\mathbb R^J\,:\,{\cal F}_{j}\geq 0,
\sum_{j=1}^{J}{\cal F}_{j}>1\right\}\;\;\mbox{and}\;\;\Xi=\left\{{\cal F}\in\mathbb R^J\,:\,{\cal F}_{j}> 0,
\sum_{j=1}^{J}{\cal F}_{j}>1\right\}.
$$

In the next theorem we describe the range of the map ${\cal F}$.
\begin{Theorem}
The following assertions hold:   {\rm (i)} ${\cal F}(\Omega)=\Xi$, {\rm (ii)} ${\cal F}(\widehat{\Omega})=\widehat{\Xi}$,\\
{\rm (iii)} ${\cal F}:\widehat{\Omega}\rightarrow \widehat{\Xi} $ is a homeomorphism.
\end{Theorem}
{\bf Proof}. First, let us prove (i). We need two estimates.

From (\ref{LUD4a}), it follows that
$$
\sum_{k=1}^J\frac{F_k(1,{\bf m})}{F_0(1,{\bf m})}=1.
$$
Therefore,
$$
\sum_{k=1}^J{\cal F}_k({\bf m})-1=\sum_{k=1}^J\frac{F_k(1,{\bf
m})}{F_0(1,{\bf
m})}\Big(\frac{1}{m_k}-1\Big)=\sum_{k=1}^J\frac{F_k(0,{\bf m}-{\bf
1})}{F_0(0,{\bf m}-{\bf 1})}\Big(\frac{1}{m_k}-1\Big).
$$
Here we have used that $F_k(1,{\bf 1})=0$ and hence $F_k(1,{\bf
m})=F_k(0,{\bf m}-{\bf 1})$, $k=0,\ldots,J$. Let  $1-m_l=\max_{1\leq k\leq J}
(1-m_k)$ and let $a_k=(1-m_k)/(1-m_l)$. In view of linearity of $F_k$ we
obtain that
$$
\sum_{k=1}^J{\cal F}_k({\bf m})-1=\sum_{k=1}^J\frac{F_k(0,{\bf
a})}{F_0(0,{\bf a})}\Big(\frac{1}{m_k}-1\Big),
$$
where ${\bf a}=(a_1,\ldots,a_J)$. We observe that $0\leq a_k\leq 1$
and $a_l=1$. Applying the maximum principle, we verify
$$
F_0(0,{\bf a})\geq F_0(0,{\bf e}_l)\neq 0
$$
and all the quantities $F_k(0,{\bf a})$ are uniformly bounded with
respect to ${\bf a}$ by linearity of $F_k$. Thus, we arrive at the
first desired inequality
\begin{equation}\label{LUD4b}
0<\sum_{k=1}^J{\cal F}_k({\bf m})-1\leq
C\sum_{k=1}^J\frac{1-m_k}{m_k}.
\end{equation}

Assuming that $\min_{1\leq k\leq J}
m_k\leq 1/2$, we prove the second inequality
\begin{equation}\label{LUD4c}
\sum_{k=1}^J{\cal F}_k({\bf m})\geq c (\min_{1\leq k\leq J}m_k)^{-1}
,
\end{equation}
where $c$ is a positive constant independent on ${\cal F}$. Let $m_l=\min_{1\leq k\leq J}m_k\leq 1/2$. Then it suffices to show that
\begin{equation}\label{LUD4d}
\frac{F_l(1,{\bf m})}{F_0(1,{\bf m})}\geq c>0.
\end{equation}
Applying the maximum principle, we get
$$
F_l(1,{\bf m})\leq F_l(1,{\bf e}_l/2).
$$
This together with boundedness of all $F_k$ yields (\ref{LUD4d}) and hence (\ref{LUD4c}).

Due to Lemmas \ref{Lem1} and \ref{Lem2}  the map
\begin{equation}\label{LUD5a}
{\cal F}\,:\,\widehat{\Omega}\to\widehat{\Xi}
\end{equation}
is open and by estimates (\ref{LUD4b}) and (\ref{LUD4c}) the map
is proper.

We fix a small positive
$\varepsilon$ and introduce
$$
\widehat{\Omega}_\varepsilon=\{{\bf m}\in
\widehat{\Omega}\,:\,m_k\geq\varepsilon,\,k=1,\ldots,J,\,\mbox{and}
\;J-m_1-\cdots-m_J\geq\varepsilon\}.
$$
Let also $\Omega_\varepsilon$ be the interior of
$\widehat{\Omega}_\varepsilon$. We represent the boundary
$\Gamma_\varepsilon$ of $\widehat{\Omega}_\varepsilon$ as
$$
\Gamma_\varepsilon=\Gamma_\varepsilon^1\cup\Gamma_\varepsilon^2\cap\Gamma_\varepsilon^3,
$$
where
$$
\Gamma_\varepsilon^1=\bigcup_{k=1}^J\{{\bf
m}\in\widehat{\Omega}_\varepsilon\,:\,F_k=0\},
$$
$$
\Gamma_\varepsilon^2=\bigcup_{k=1}^J\{{\bf
m}\in\widehat{\Omega}_\varepsilon\,:\,m_k=\varepsilon\}
$$
and
$$
\Gamma_\varepsilon^3=\left\{{\bf
m}\in\widehat{\Omega}_\varepsilon\,:\,J-\sum_{k=1}^Jm_k=\varepsilon\right\}.
$$
Then ${\cal F}(\Gamma_\varepsilon^1)\subset \partial\Xi$ and
due to estimates (\ref{LUD4b}) and (\ref{LUD4c})
$$
{\cal F}(\Gamma_\varepsilon^2)\subset\left\{{\cal
F}\in\widehat{\Xi}\,:\,|{\cal F}|\geq  c\varepsilon^{-1}\right\}
$$
and
$$
{\cal F}(\Gamma_\varepsilon^3)\subset\left\{{\cal
F}\in\widehat{\Xi}\,:\,\sum_{k=1}^J{\cal F}_k-1\leq c\varepsilon\right\},
$$
where $c$ is a positive constant independent of $\varepsilon$.
Consider connected components of $\mathbb R^J\setminus {\cal
F}(\Gamma_\varepsilon)$. Due to the above-performed analysis of ${\cal
F}(\Gamma_\varepsilon^j)$, one of these components contains
the set
\begin{equation}\label{LUDbg}
\left\{{\cal F}\in\Xi\,: \,|{\cal F}|<
\frac{c}{\e},\;\sum^{J}_{j=1}{\cal F}_{j}-1> c\e \right\}.
\end{equation}
Since the map is open, each connected component must belong to
${\cal F}(\Omega_\varepsilon)$ or its intersection with the last
set is empty. Thus we conclude that the set (\ref{LUDbg}) lies in ${\cal
F}(\Omega_\varepsilon)$. Sending $\varepsilon$ to $0$, we obtain
$$
{\cal F}(\Omega)=\Xi,
$$
which proves (i).

Since the proper map ${\cal F}$ is a local homeomorphism, we
arrive at (ii).

Assertion (iii) follows from assertion (ii) and Lemma  \ref{Lem2}. The proof is completed.

\subsection{The case of permeable walls, $B\geq 0$}\label{Sec.3.2}

In this section we assume that  $B\geq 0$ and $B$ is not identically zero. This guarantees that the form $a(u,u)$ is positive definite on $H^1(G^0)$. This allows us to introduce the function $Q$ as the solution to equation (\ref{00.2}) with boundary conditions
\begin{equation}\label{K1b}
Q(W_0)=1\;\; \mbox{and}\;\; F_j(Q)=0\,\; \mbox{for}\;\; j=1,\ldots,J.
\end{equation}
Repeating the proof of Lemma \ref{QuQu1}, we verify that this solution is positive.

In our study of equation (\ref{03newQ})  an important role will be played by the following

\smallskip
\noindent {\bf Maximum principle  ($B\geq 0$)}. {\em Let $v$ be a
solution to the problem (\ref{00.2})-(\ref{00.3}) with the Dirichlet
data $(X_{0},{\bf X})$ and let $v$ be not identically constant. If $\max_{G^0}v\geq 0$, then
$$
 v< \max_{0\leq j\leq
J}X_{j}\;\;\;\mbox{in the interior of $G^0$.}
$$
Moreover,  if the maximum is attained  at $W_j$ then $F_j<0$
when $j=1,\ldots,J$ and $F_j>0$ when $j=0$.

Similar assertion is valid for the minimum. (It suffices to apply the above principle to the function $-v$).}

\medskip
This principle is well known for second order elliptic partial differential equations, see \cite{PW}, Th.6 and 8. The graph version is also well known and its proof is quite straightforward.

\medskip
Let us examine the function $Q$. Since $Q$ is a positive function, the application of the maximum principle gives
\begin{equation}\label{K3a}
Q(W_k)<1,\;\;k=1,\ldots,J, \;\;\mbox{and}\;\; F_0(Q)>0.
\end{equation}
Next, we show that the sets (\ref{K3}) and (\ref{K2}) can be described as
\begin{equation}\label{XX2a}
\Omega=\{ {\bf m}\;|\;m_j>0\;\mbox{and}\;F_{j}(1,{\bf m})> 0,\; j=1,\ldots, J\}
\end{equation}
and
\begin{equation}\label{XX2b}
\widehat{\Omega}=\{ {\bf m}\;|\;m_j>0\;\mbox{and}\;F_{j}(1,{\bf m})\geq 0,\; j=1,\ldots, J\}.
\end{equation}
Indeed, let ${\bf m}$ belong to the right-hand side of (\ref{XX2a}). Applying the maximum principle to the function $M-Q$, where $M $ is the solution to (\ref{00.2}) with the Dirichlet data $(1,{\bf m})$, we conclude that $M-Q$ is negative in $\overline{G^0}\setminus \{W_0\}$ (if $\max_{1\leq j\leq J}(M-Q)(W_j)\geq 0$ then the flux must be negative at the point, where this maximum is attained, but $F_j(M-Q)>0$  and hence
\begin{equation}\label{K3b}
m_j< Q(W_j),\;\;j=1,\ldots,J,\;\;\; F_0(1,{\bf m})> F_0(Q).
\end{equation}
Relation (\ref{XX2b}) is proved similarly.

\begin{Lemma}\label{Lem2a} The Jacobian ${\cal J}({\bf m})$ is invertible at any point ${\bf m} \in \widehat{\Omega}$. The
 map ${\cal F}:\widehat{\Omega}\rightarrow \mathbb{R}^{J}$ is injective.
\end{Lemma}
{\bf Proof}. The proof is literary the same as that of Lemma \ref{Lem2}.

\medskip
Now we are in position to describe all solutions to (\ref{03newQ}).

\begin{Theorem}  The following assertion hold:

\noindent
  {\rm (i)} ${\cal F}(\Omega)=\mathbb R_+^J=\{{\cal F}\in\mathbb R^J\,:\, {\cal F}_j>0\}$,\\
   {\rm (ii)} ${\cal F}(\widehat{\Omega})=\overline{\mathbb R_+^J}$,\\
{\rm (iii)} The map
\begin{equation}\label{KK2s}
{\cal F}:\widehat{\Omega}\mapsto \overline{\mathbb R_+^J}
\end{equation}
 is a homeomorphism.
\end{Theorem}
{\bf{Proof}}.
Let us prove (i).
Due to Lemma \ref{Lem2a}  the map $\Omega\ni{\bf m}\mapsto{\cal F}({\bf m})$
is open. We fix a small positive
$\varepsilon$ and introduce
$$
\widehat{\Omega}_\varepsilon=\{{\bf m}\in
\widehat{\Omega}\,:\,m_k\geq\varepsilon F_k(1,{\bf m}),\,k=1,\ldots,J\}.
$$
By maximum principle $m_k=0$ for a certain $k$ implies $F_k(1,{\bf m})\neq 0$. Therefore $\widehat{\Omega}_\varepsilon$ is compact in $\widehat{\Omega}$. Let  $\Omega_\varepsilon$ be the interior of
$\widehat{\Omega}_\varepsilon$. We represent the boundary
$\Gamma_\varepsilon$ of $\widehat{\Omega}_\varepsilon$ as
$$
\Gamma_\varepsilon=\Gamma_\varepsilon^1\cup\Gamma_\varepsilon^2,
$$
where
$$
\Gamma_\varepsilon^1=\bigcup_{k=1}^J\{{\bf
m}\in\widehat{\Omega}_\varepsilon\,:\,F_k(1,{\bf m})=0\},\;\;\;
\Gamma_\varepsilon^2=\bigcup_{k=1}^J\{{\bf
m}\in\widehat{\Omega}_\varepsilon\,:\,m_k=\varepsilon F_k(1,{\bf m})\}.
$$
Then ${\cal F}(\Gamma_\varepsilon^1)\subset \bigcup_{k=1}^J\{{\cal F}\in \mathbb R^J\,:\, {\cal F}_k=0\}$ and
due to the second  estimate in (\ref{K3b})
$$
{\cal F}(\Gamma_\varepsilon^2)\subset\{{\cal
F}\in \overline{\mathbb R_+^J}\,:\,|{\cal F}|\geq  c\varepsilon^{-1}\}.
$$
Consider connected components of $\mathbb R^J\setminus {\cal
F}(\Gamma_\varepsilon)$. Due to the above-perfomed analysis of ${\cal
F}(\Gamma_\varepsilon^j)$, one of connected components contains
the set
\begin{equation}\label{LUDb}
\left\{{\cal F}\in\mathbb R_+^J\,: \,|{\cal F}|<
\frac{c}{\e} \right\}.
\end{equation}
Since the map is open, each of connected components must belong to
${\cal F}(\Omega_\varepsilon)$ or its intersection with this
set is empty. Thus we get that the set (\ref{LUDb}) lies in ${\cal
F}(\Omega_\varepsilon)$. Sending  $\varepsilon$ to $0$, we obtain
$$
{\cal F}(\Omega)=\mathbb R_+^J,
$$
that proves (i).

Since the proper map ${\cal F}$ is a local homeomorphism, we
arrive at (ii).

Assertion (iii) follows from assertion (ii) and Lemma  \ref{Lem2a}. The proof is completed.

\subsection{The case of permeable walls, $B\leq 0$}\label{Sec.3.3}

In this section we consider the case  $B\leq0$ and $B$ is not identically zero. 

We denote by $R$ the solution to the Dirichlet problem (\ref{00.2}), (\ref{00.3}) with $X_0=\cdots=X_J=1$. By Lemma \ref{QuQu1} this solution is positive.

We will use the following sharpening of Lemma \ref{QuQu1}.
\begin{Lemma}\label{QuQu1a} Let $w\in H^1(G^0)$ satisfy
\begin{equation}\label{00.2a}
a(w,u)=\sum_{e\in E}\int_0^{L_e}F_e(\xi)\,u(\xi)d\xi.\;\;\;\mbox{for all $u\in H^1_0(G^0)$}
\end{equation}
and {\rm (\ref{00.3})}, where $X_j\geq 0$, $j=0,\ldots,J$, and the functions $F_e$ are continuous and bounded on $e$. If $w$ is not identically zero, then
 $w>0$ in $G^0$. If additionally $X_k=0$ for certain $k=0,\ldots,J$, then $F_k(X_0,{\bf X})<0$ for $k=0$ and $F_k(X_0,{\bf X})>0$ for $k>0$.
\end{Lemma}
The proof repeats the proof of Lemma \ref{QuQu1}. But since we cannot use the uniqueness for the Cauchy problem in this case one should use instead the following property:
if a non-negative function, satisfying (\ref{00.2a}), is zero at a certain point, then its derivative is non-negative before this point and non-positive after this point.

\medskip
Applying this lemma to the function $R$ we obtain that
\begin{equation}\label{PW2}
R>1,\;\;\mbox{in $G^0$}\;\; F_0(1,{\bf 1})<0\;\;\;\mbox{and}\;\;  F_j(1,{\bf 1})>0,\;\; j=1,\ldots,J.
\end{equation}

We assume in this section that $F_0(1,{\bf 0})>0$. Since $F_j(1,{\bf 0})>0$, $j=1,\ldots,J$, we see that small ${\bf m}$ with positive components belong to $\Omega$.

In the next lemma we study local invertibility of the Jacobian.

\begin{Lemma} Let {\rm (\ref{TTf1s})} be valid for any ${\bf m}\in\widehat{\Omega}$ and non vanishing ${\bf H}\in{\cal R}_0$.
Then the Jacobian ${\cal J}(\bf m)$ is not degenerate at any point ${\bf m} \in \widehat{\Omega}$
 outside the surface ${\cal S}$:
 \begin{equation}\label{PW3}
 S({\bf m}):=F_0(1,{\bf m})-\sum_{j=1}^{J}F_{j}(1,{\bf m})m_{j}=0.
 \end{equation}
 This surface $S$ satisfies
  \begin{equation}\label{PW4}
 \sum_{j=1}^J m_j\partial_{m_j}S({\bf m})<0\;\;\mbox{when $S=0$},
 \end{equation}
 which means that ${\cal S}$ is star-shaped with respect to the origin.
\end{Lemma}
{\bf Proof.} The first assertion follows from Lemma \ref{Lem1}.

To prove the second assertion (inequality (\ref{PW4}))  we note that
$$
\sum_{j=1}^J m_j\partial_{m_j}S({\bf m})=F_0(0,{\bf m})-\sum_{j=1}^J(F_j(1,{\bf m})m_j+F_j(0,{\bf m})m_j).
$$
Using that $S({\bf m})=0$ we get
\begin{eqnarray*}
&&\sum_{j=1}^J m_j\partial_{m_j}S({\bf m})=-F_0(1,{\bf 0})-\sum_{j=1}^JF_j(0,{\bf m})m_j\\
&&=-F_0(1,{\bf 0})-F_0(1,{\bf m})+\sum_{j=1}^JF_j(1,{\bf 0})m_j=-2F_0(1,{\bf m}).
\end{eqnarray*}
Since $F_0(1,{\bf m})$ is positive the proof is complete.

We introduce
$$
\widehat{\Omega}_+=\{{\bf m}\in\widehat{\Omega}\,:\, F_0(1,{\bf m})-\sum_{j=1}^{J}F_{j}(1,{\bf m})m_{j}>0\}.
$$
\begin{Theorem}\label{Omega+} Let {\rm (\ref{TTf1s})} be valid for any ${\bf m}\in\widehat{\Omega}$ and non vanishing ${\bf H}\in{\cal R}_0$.
 Then the map ${\cal F}$ is injective on $\widehat{\Omega}_+ $. The map is not injective on $\Omega $.
\end{Theorem}
{\bf Proof.} The first assertion follows from Lemma \ref{Lem1s}.

To prove the second assertion we start from the  identity
\begin{equation}\label{Dec24a}
{\cal F}_j({\bf m})-{\cal F}_j({\bf n})=\frac{F_j(0,{\bf h})-F_0(1,{\bf n}){\cal F}_j({\bf n})h_j-{\cal F}_j({\bf n})m_jF_0(0,{\bf h})}{F_0(1,{\bf m})m_j},
\end{equation}
where ${\bf h}={\bf m}-{\bf n}$. We represent ${\bf h}$ as ${\bf h}={\bf H}+\lambda{\bf n}$ where ${\bf H}\in{\cal R}_0$ and $\lambda\in\mathbb R$. Then
 the equalities ${\cal F}_j({\bf m})={\cal F}_j({\bf n})$, $j=1,\ldots,J$, are equivalent to
\begin{equation}\label{Dec29a}
F_j(0,{\bf H})-F_0(1,{\bf n}){\cal F}_j({\bf n})H_j=\lambda\big(F_j(1,{\bf 0})+{\cal F}_j({\bf n})m_jF_0(0,{\bf n})\big).
\end{equation}
where we have used the identity $F_0(1,{\bf n}){\cal F}_j({\bf n})n_j=F_j(1,{\bf n})$. One solution to
 (\ref{Dec29a}) is $\lambda=0$ and ${\bf H}=0$. Let us find another solution. Multiplying both sides of (\ref{Dec29a}) by $n_j$ and summing up, we get
$$
\sum_{j=1}^Jn_j\big(F_j(1,{\bf 0})+{\cal F}_j({\bf n})m_jF_0(0,{\bf n})\big)=0.
$$
Since $\sum n_jF_j(1,{\bf 0})=-F_0(0,{\bf n})$, the last relation takes the form
$$
F_0(0,{\bf n})\Big(1-\sum_{j=1}^J{\cal F}_j({\bf n})m_jn_j\Big)=0,
$$
which leads to
\begin{equation}\label{Dec29b}
\lambda=\rho-\frac{\sum_{j=1}^J{\cal F}_j({\bf n})n_jH_j}{\sum_{j=1}^J{\cal F}_j({\bf n})n_j^2},\;\;\rho=\frac{1-\sum_{j=1}^J{\cal F}_j({\bf n})n_j^2}{\sum_{j=1}^J{\cal F}_j({\bf n})n_j^2}.
\end{equation}
Now we can write equation (\ref{Dec29a}) as
\begin{eqnarray}\label{Jan2a}
&&-L_j({\bf n};{\bf H})+\frac{\sum_{k=1}^J{\cal F}_k({\bf n})n_kH_k}{\sum_{k=1}^J{\cal F}_k({\bf n})n_k^2}\cdot{\cal F}_j({\bf n})h_jF_0(0,{\bf n})\nonumber\\
&&=\rho\big(F_j(1,{\bf 0})+{\cal F}_j({\bf n})m_jF_0(0,{\bf n})\big),
\end{eqnarray}
where
\begin{eqnarray*}
&&L_j({\bf n};{\bf H})=-F_j(0,{\bf H})+F_0(1,{\bf n}){\cal F}_j({\bf n})H_j\\
&&-\frac{\sum_{k=1}^J{\cal F}_k({\bf n})n_kH_k}{\sum_{k=1}^J{\cal F}_k({\bf n})n_k^2}\cdot (F_j(1,{\bf 0})+{\cal F}_j({\bf n})n_jF_0(0,{\bf n})\big).
\end{eqnarray*}
Further, we note that
\begin{eqnarray}\label{Jan5a}
&&\sum_{j=1}^JL_j({\bf n};{\bf H})H_j=\sum_{j=1}^J\big(-F_j(0,{\bf H})H_j+F_0(1,{\bf n}){\cal F}_j({\bf n})H_j^2\big)\nonumber\\
&&-\frac{F_0(0,{\bf n})}{\sum_{k=1}^J{\cal F}_k({\bf n})n_k^2}\big(\sum_{k=1}^J{\cal F}_k({\bf n})n_kH_k\big)^2.
\end{eqnarray}
Since $F_0(0,{\bf n})<0$, the above quadratic form is positive for nontrivial ${\bf H}$ due to (\ref{TTf1s}). Therefore we can
apply the fixed point theorem to obtain an existence of solution ${\bf H}$ to system (\ref{Jan2a}) satisfying ${\bf H}=O(\rho)$. The same relation is valid
for $\lambda$ because of (\ref{Dec29b}). To show that $\lambda$ is not zero we derive from (\ref{Jan5a}) the existence of
 a positive constant $\alpha=\alpha(n)$ such that
\begin{equation}\label{Jan2b}
\sum_{j=1}^JL_j({\bf n};{\bf H})H_j\geq \frac{(1+\alpha)|F_0(0,{\bf n})|}{\sum_{k=1}^J{\cal F}_k({\bf n})n_k^2}X^2,\;\;X=\sum_{k=1}^J{\cal F}_k({\bf n})n_kH_k
\end{equation}
for ${\bf H}\in {\cal R}_0\setminus{\cal O}$ and ${\bf n}\in\widehat{\Omega}$. Multiplying both sides of (\ref{Jan2a}) by $H_j$ and summing up and using (\ref{Jan2b}), we get
$$
\frac{(1+\alpha)|F_0(0,{\bf n})|}{\sum_{k=1}^J{\cal F}_k({\bf n})n_k^2}X^2\leq \rho(|F_0(0,{\bf n})|\,|X|+\sum_{j=1}^J{\cal F}_j({\bf n})h_jH_jF_0(0,{\bf n}))
$$
which leads to the estimate
$$
\frac{|X|}{\sum_{k=1}^J{\cal F}_k({\bf n})n_k^2}\leq \frac{\rho+C\rho^2}{(1+\alpha)}
$$
where $C$ is a constant depending on ${\bf n}$. Using this inequality together with (\ref{Dec29b}), we conclude that
$$
\lambda\geq\frac{\alpha\rho}{1+\alpha}+O(\rho^2).
$$
Thus the pair $\lambda,\,{\bf H}$ give another non-zero solution to ${\cal F}({\bf m})={\cal F}({\bf n})$. The proof is complete.

\subsubsection{An example}

Here we present an example of non uniqueness of solution to problem (\ref{00.2}), (\ref{00.3}) in the case of $B\leq0$.

Consider one-dimensional case, i.e. the graph $G^0$ has $J=1$ output. Let $B=-\gamma^2\leq0$ ($\gamma$ is not identically zero) and $H\equiv1$ in (\ref{strong}). We look for the solution $w(\xi)$ to problem (\ref{00.2}), (\ref{00.3}) on the interval $[0,1]$ with the Dirichlet data $(1,m)$. This solution is given by
$$
w(\xi)=\cos\gamma\xi-A\sin\gamma\xi,\;\;\;A=\frac{\cos\gamma-m}{\sin\gamma},
$$
where $m\in (0,\cos\gamma)$ and $0<\gamma <\pi$. Direct calculations lead to
$$
{\cal F}(m)=\frac{1-m {\rm cos}\gamma}{({\rm cos}\gamma-m)m},
$$
It is easy to see that the map ${\cal F}(m): (0,{\rm cos}\gamma)\rightarrow(0,+\infty)$ is not a homeomorphism.

\begin{figure}[ht!]
\center{\includegraphics[scale=0.45]{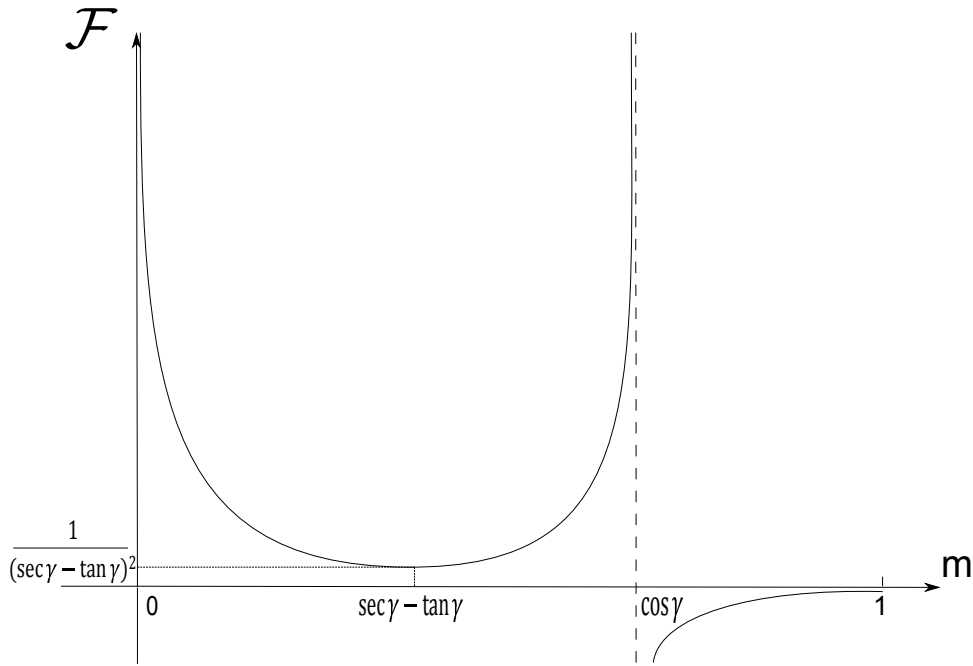}}\\ \caption{Dependence of ${\cal F}$
on $m$ (the case of $J=1$ output)}\label{1D}
\end{figure}

If $m \rightarrow 0$ we obtain
$$
{\cal F}=\frac{1}{m{\rm cos}\gamma }+O(1).
$$
If $m \rightarrow {\rm cos}\gamma$ we have
$$
{\cal F}=\frac{1-{\rm cos}^2\gamma}{({\rm cos}\gamma-m){\rm cos}\gamma}+O(1).
$$
This shows that the equation
$${\cal F}(m)=\kappa$$
has  two solutions on $(0,{\rm cos}\gamma)$ for $\kappa>\min_{0<m<{\rm cos}\gamma}{\cal F}(m)$.


\section{Appendix}
\subsection{The Reynolds equation for a viscous incompressible flow}
In a thin tube with curved walls
\begin{equation}\label{A1}
\Omega^{h}=\left\{x=(y,z)\in \mathbb{R}^{2}\times \mathbb{R}^{1}: \eta:=h^{-1}y \in \omega(z), z\in (0,1)\right\}
\end{equation}
we consider the Navier-Stokes equations
\begin{equation}\label{A2}
-\Delta_{x}u^{h}(x)+{\sf Re}(u^{h}(x)\cdot\nabla_{x})u^{h}(x)+\nabla_{x}p^{h}(x)=0,\quad x\in\Omega^{h},
\end{equation}
\begin{equation}\label{A3}
-\nabla_{x}\cdot u^{h}(x)=0,\quad x \in \Omega^{h},
\end{equation}
where $u^{h}$ and $p^{h}$ are a velocity vector and a pressure respectively, $\nabla_{x}=\mbox{grad}, \nabla_{x}\cdot=\mbox{div}, \Delta_{x}=\nabla_{x}\cdot\nabla_{x}$ is the Laplace operator. Notice that, by rescaling, length of the tube has been reduced to 1, equations (\ref{A2}), (\ref{A3}) are written in the dimensionless form and involve the Reynolds number {\sf Re} which is compared with the small parameter $h \in (0,1]$ as follows:
\begin{equation}\label{A4}
{\sf Re}=h^{\rho}R_{0}.
\end{equation}
Furthermore, $\omega(z)=\k_{z}\omega$ where $\omega$ is a domain in the plane $\mathbb{R}^{2}\ni\eta=(\eta_{1},\eta_{2})$ bounded by a simple smooth closed contour $\partial \omega$ and $\k_{z}: \mathbb{R}^{2}\rightarrow\mathbb{R}^{2}$ is a family of diffeomorphisms dependent smoothly on the longitudinal coordinate $z=x_{3}\in [0,1]$.

At the lateral boundary $\Sigma^{h}=\left\{x:\eta\in \partial\omega(z), z\in (0,1)\right\}$, we impose the no-slip condition
\begin{equation}\label{A5}
u^{h}_{s}(x)=0,\quad x\in\Sigma^{h},
\end{equation}
and the percolation condition
\begin{equation}\label{A6}
u^{h}_{n}(x)=\beta^{h}(x)(p^{h}(x)-\e_{nn}(u^{h};x)),\quad  x\in \Sigma^{h},
\end{equation}
where $u^{h}_{n}=n^{h}\cdot u^{h}, n^{h}$ is the unit vector of the outward normal on $\Sigma^{h}, u^{h}_{s}=u^{h}-u^{h}_{n}n^{h}$ is a two-dimensional vector of tangential velocities, and $\e=(\e_{jk})$ is the shear stress tensor with components
$$
\e_{jk}(u^{h})=\frac{1}{2}\left(\frac{\partial u^{h}_{j}}{\partial x_{k}}+\frac{\partial u^{h}_{k}}{\partial x_{j}}\right),\  j,k=1,2,3,\  \e_{nn}(u^{h})=n^{h}\cdot \e(u^{h})n^{h}.
$$

The boundary condition (\ref{A6}) means that a percolation occurs through the wall $\Sigma^{h}$ and is proportional to the normal hydrodynamic force with the coefficient (positive or negative)
\begin{equation}\label{A7}
\beta^{h}(x)=h^{3}b(\eta,z).
\end{equation}

We are looking for an asymptotic solution of problem (\ref{A2}), (\ref{A3}), (\ref{A5}), (\ref{A6}) in the form
\begin{equation}\label{A8}
p^{h}(x)=h^{-4}P(z)+h^{-3}Q(\eta,z)+\ldots,
\end{equation}
\begin{equation}\label{A9}
u^{h}(x)=(h^{-1}V_{1}(\eta,z),h^{-1}V_{2}(\eta,z),h^{-2}W(\eta,z))+\ldots.
\end{equation}
For the normal at $\Sigma^{h}$ we have the formula
\begin{equation}\label{A10}
n^{h}(x)=(1+h^{2}|N_{0}(\eta,z)|^2)^{-\frac{1}{2}}(N_{1}(\eta,z),N_{2}(\eta,z),hN_{0}(\eta,z)),
\end{equation}
where $N=(N_{1},N_{2})$ is the unit outward normal on the boundary of the domain $\omega(z)\subset \mathbb{R}^{2}$ while the component $hN_{0}$ reflects the variability of the tube cross-section $\omega^{h}(z),$ see (\ref{A13}) below.

Assuming the Reynolds number to be small that is $\rho>0$ in (\ref{A4}), we insert (\ref{A7})-(\ref{A10}) into (\ref{A2}),(\ref{A3}) and (\ref{A5}),(\ref{A6}). Then we collect coefficients of like powers of $h$ and compose two problems
\begin{equation}\label{A11}
\left.\begin{array}c -\Delta_{\eta}W(\eta,z)=-\partial_{z}P(z), \eta \in \omega(z),\\
W(\eta,z)=0, \eta \in \partial \omega(z),
\end{array}
\right.
\end{equation}
and
\begin{equation}\label{A12}
\left.\begin{array}c -\Delta_{\eta}V(\eta,z)+\nabla_{\eta}Q(\eta,z)=0, \eta \in \omega(z),\\
-\nabla_{\eta}V(\eta,z)=\partial_{z}W(\eta,z), \eta \in \omega(z),\\
V_{N}(\eta,z)=b(\eta,z)P(z)-N_{0}(\eta,z)W(\eta,z), \eta \in \partial \omega(z),\\
V_{S}(\eta,z)=0, \eta \in \partial \omega(z).
\end{array}
\right.
\end{equation}
The solution of (\ref{A11}) is
\begin{equation}\label{A13}
W(\eta,z)=-\frac{1}{2}\Psi(\eta,z)\partial_{z}P(z),
\end{equation}
where $\Psi$ is the Prandtl function satisfying
$$
-\Delta_{\eta}\Psi(\eta,z)=2, \eta \in \omega(z),\quad \Psi(\eta,z)=0, \eta \in \partial \omega(z).
$$
Moreover, the formula
\begin{equation}\label{A14}
\frac{d}{dz}\int_{\omega(z)}\Phi(\eta,z)d\eta=\int_{\omega(z)}\frac{\partial \Phi}{\partial z}(\eta,z)d\eta-\int_{\partial \omega(z)}N_{0}(\eta,z)\Phi(\eta,z)dS_{\eta}
\end{equation}
for differentiation of integrals with variable limits transforms a compatibility condition (the total flux vanishes) in the two-dimensional Stokes problem (\ref{A12}) into the ordinary differential equation of Reynolds type
\begin{equation}\label{A15}
-\partial_{z}\left(H(z)\partial_{z}P(z)\right)+B(z)P(z)=0, z \in (0,1)
\end{equation}
(see, e.g., \cite{139} for details). Here, $4H(z)$ is the torsion rigidity of the domain $\omega(z)$, see, e.g., \cite{PoSe}, and $B(z)$ stands for the total percolation coefficient, namely
\begin{equation}\label{A16}
H(z)=\frac{1}{2}\int_{\omega(z)}\Psi(\eta,z)d\eta=\frac{1}{4}\int_{\omega(z)}|\nabla_{\eta}\Psi(\eta,z)|^{2}d\eta, B(z)=\int_{\partial \omega(z)}b(\eta,z)dS_{\eta}.
\end{equation}

According to (\ref{A9}) and (\ref{A13}), (\ref{A16}), the flux through the cross-section $\omega^{h}(z)$ of the tube is determined as
\begin{equation}\label{A17}
F^{h}(z)=\int_{\omega^{h}(z)}u^{h}_{3}(y,z)dy=h^{2}\int_{\omega(z)}h^{-2}W(\eta,z)d\eta+\ldots=-H(z)\partial_{z}P(z)+\ldots
\end{equation}
In this way, the thin tube (\ref{A1}) is able to drive flux of order $1=h^{0}$ within the linear one-dimensional model under restriction (\ref{A4}) with $\rho>0$ of the small Reynolds number. However, if the flux is infinitesimal as $h\rightarrow +0$, the restriction on {\sf Re} can be weakened. These observations, owing to (\ref{A9}) and (\ref{A4}), are supported by calculation of the convective term in (\ref{A2})
$$
{\sf Re}(u^{h}(x)\cdot\nabla_{x})u^{h}(x)=R_{0}h^{\rho-4}(V(\eta,z)\cdot\nabla_{\eta}+W(\eta,z)\partial_{z})(0,0,W(\eta,z))+\ldots
$$
At $\rho=0$, the latter term must come to problem (\ref{A11}) that deprives the first limit problem of sense. In other words, to validate a linear one-dimensional model, one has either to assume $\rho>0$ in (\ref{A4}), or to reduce negative exponents of the small parameter $h$ in the asymptotic ansatzes 
(\ref{A9}) and (\ref{A8}). We also mention that an intensive enforced percolation may change asymptotic structures of a thin flow, cf. \cite{334}.

\subsection{The one-dimensional flow of an ideal liquid}

Let now $u^h (x)$ be the velocity potential in an ideal liquid which satisfies the Laplace equation in thin curved tube (\ref{A1})
\begin{equation}\label{A21}
-\Delta_{x}u^h(x)=0,\quad x \in \Omega^h,
\end{equation}
as well as the boundary condition of Robin's type with the coefficient (\ref{A7}) on the lateral surface
\begin{equation}\label{A22}
\partial_{n^h}u^h(x)+\beta^h(x)u^h(x)=0,\quad x \in \Sigma^h.
\end{equation}
The latter is intendent to describe the percolation law through the wall. We assume the standard asymptotic ansatz, cf. \cite{MaNaPl},  Ch.15,16,
\begin{equation}\label{A23}
u^h(x)=h^{-2}U(z)+W(\eta,z)+\ldots
\end{equation}

Inserting (\ref{A23}), (\ref{A7}) into (\ref{A21}), (\ref{A22}) and extracting coefficient of $h^{-2},h^{-1}$ respectively yield the planar Neumann problem on the cross-section

\begin{equation}\label{A24}
\left.\begin{array}c -\Delta_\eta W(\eta,z)=\partial^2_z U(z), \eta \in \omega(z),\\
\partial_N W(\eta,z)=-N_0(\eta,z)\partial_z U(z)-b(\eta,z)U(z), \eta \in \partial \omega(z).
\end{array}
\right.
\end{equation}
The compatibility condition in the problem (\ref{A24}) reads:
\begin{equation}\label{A25}
\int_{\omega(z)}\partial^2_z U(z)d\eta-\int_{\partial \omega(z)}N_0(\eta,z)dS_\eta \partial_z U(z)-\int_{\partial \omega(z)}b(\eta,z)dS_\eta U(z)=0.
\end{equation}
Formula (\ref{A14}) with $\Phi=1$ becomes
\begin{equation}\label{A26}
-\int_{\partial \omega(z)}N_0(\eta,z)dS_\eta =\frac{d}{dz}\int_{\omega(z)}d\eta=\frac{dH}{dz}(z),
\end{equation}
where $H(z)=|\omega(z)|$ stays for area of the figure $\omega(z)\subset \mathbb{R}^2$. Using (\ref{A26}) and the second definition in (\ref{A16}) converts (\ref{A25}) into the ordinary differential equation
\begin{equation}\label{A27}
-\partial_z\left(H(z)\partial_z U(z)\right)+B(z)U(z)=0,\quad z \in (0,1).
\end{equation}

According to (\ref{A23}) the projection onto the z-axis of the velocity vector $\nabla u^h(x)$ is $\partial_z u^h(x)=h^{-2}\partial_z w(z)$ so that the flux through the cross-section $\omega^h(z)$ of the tube (\ref{A1}) is equal to
\begin{equation}\label{A28}
\int_{\omega^h (z)}\partial_z u^h(y,z)dy=H(z)\partial_z U(z)+\ldots
\end{equation}

\bigskip
\noindent {\bf Acknowledgements.} V.~K. acknowledges
the support of the Swedish Research Council (VR) grant EO418401. S.~N. was supported by the Russian Foundation for Basic Research, project no.12-01-00348, and by
  Link\"oping University (Sweden).   G.~Z. was supported by  Link\"oping University and RFBR grants 15-31-20600, 16-31-60112.


\begin{thebibliography}{20}

\bibitem{BCG}
{\it Bressan, Alberto; Canic, Suncica; Garavello, Mauro; Herty, Michael; Piccoli, Benedetto},
Flows on networks: recent results and perspectives. (English summary)
EMS Surv. Math. Sci.  1  (2014),  no. 1, 47–111.


\bibitem{APQ} {\it D'Angelo, Carlo; Panasenko, Grigory; Quarteroni, Alfio}, Asymptotic-numerical derivation of the Robin type coupling conditions for the macroscopic pressure at a reservoir-capillaries interface, Appl. Anal.  92  (2013),   1, 158–171.

    \bibitem{Dav} Jeffrey M. Davis, On the Linear Stability of Blood Flow Through Model Capillary Networks, Bulletin of Mathematical Biology, December 2014, Volume 76, Issue 12,  pp 2985–3015.


\bibitem{FiLaLi} G. Fibich, Y. Lanir, N. Liron, Mathematical model of blood flow in a coronary capillary,
American Journal of Physiology - Heart and Circulatory Physiology, 1993  Vol. 265  no.  5,   H1829-H1840.

\bibitem{PW} Protter, Murray H.; Weinberger, Hans F. Maximum principles in differential equations. Corrected reprint of the 1967 original. Springer-Verlag, New York, 1984.

\medskip

\bibitem{139} {\it Nazarov S.A., Pileckas K.I.}, Reynolds flow of a
fluid in a thin three-dimensional channel, Litovsk. mat. sbornik.
1990. V. 30, N 4. P. 772-783  (English transl.: Lithuanian Math. J.
1990. V. 30, N 4. P. 366-375).

\medskip

\bibitem{334} {\it Nazarov S.A., Videman J.H.}, Reynolds type equation
for a thin flow under intensive transverse percolation, Math.
Nachr. 2004. Bd. 269/270. S. 189-209.

\medskip

\bibitem{MaNaPl} {\it Maz'ya V., Nazarov S., Plamenevskij B.} Asymptotic
theory of elliptic boundary value problems in singularly perturbed
domains. Vol. 2. Basel: Birkh\"auser Verlag, 2000.

\medskip

\bibitem{PoSe} {\it Polya, G.; Szeg\"o, G.} Isoperimetric Inequalities in Mathematical Physics. Annals of Mathematics Studies, no. 27, Princeton University Press, Princeton, N. J., 1951.

\bibitem{Poz} C. Pozrikidis, Numerical Simulation of Blood Flow Through Microvascular Capillary Networks, Bulletin of Mathematical Biology, August 2009, Volume 71, Issue 6,  pp 1520–1541.




\end{thebibliography}
\end{document}